\documentclass[11pt]{article}
\RequirePackage{amsthm,amsmath,amsfonts,amssymb}
\RequirePackage[colorlinks,linkcolor=red,citecolor=blue,urlcolor=blue]{hyperref} 
\usepackage[authoryear]{natbib}
\usepackage{amsfonts}
\usepackage{bbm}
\usepackage{mathrsfs}
\usepackage{mathtools}
\usepackage{subfigure}
\usepackage{titlesec}
\usepackage{lipsum}
\usepackage{enumerate}
\usepackage{color}
\usepackage{amssymb}
\usepackage{amsmath,amssymb,amsfonts,amsthm}
\usepackage{bm}
\usepackage{here}
\usepackage{multicol}
\usepackage{epsf}
\usepackage{authblk}
\usepackage{tikz}
\usetikzlibrary{positioning}
\usetikzlibrary{arrows.meta}
\makeatother
\begin{document}
\newcommand{\field}[1]{\mathbb{#1}} 
\newcommand{\ind}{\mathbbm{1}}
\newcommand{\C}{\mathbb{C}}
\newcommand{\D}{\,\mathscr{D}}
\newcommand{\E}{\,\mathrm{E}}
\newcommand{\Prob}{\,\mathrm{P}}
\newcommand{\F}{\,\mathscr{F}}
\newcommand{\I}{\,\mathrm{i}}
\newcommand{\N}{\field{N}}
\newcommand{\T}{\,\mathrm{T}}
\newcommand{\Ls}{\,\mathscr{L}}
\newcommand{\ve}{\,\mathrm{vec}}
\newcommand{\var}{\,\mathrm{var}}
\newcommand{\cov}{\,\mathrm{Cov}}
\newcommand{\vech}{\,\mathrm{vech}}
\newcommand*\dif{\mathop{}\!\mathrm{d}}
\newcommand{\difs}{\mathrm{d}}
\newcommand\mi{\mathrm{i}}
\newcommand\me{\mathrm{e}}
\newcommand{\R}{\field{R}}
\newcommand{\es}{\hat{f}_n}
\newcommand{\ess}{\hat{f}_{\flr{ns}}}
\newcommand{\p}[1]{\frac{\partial}{\partial#1}}
\newcommand{\pp}[1]{\frac{\partial^2}{\partial#1\partial#1^{\top}}}
\newcommand{\para}{\bm{\theta}}
\providecommand{\bv}{\mathbb{V}}
\providecommand{\bu}{\mathbb{U}}
\providecommand{\bt}{\mathbb{T}}
\newcommand{\flr}[1]{\lfloor#1\rfloor}
\newcommand{\ba}{B_m}
\newcommand{\bxi}{\bar{\xi}_n}
\newcommand{\sgn}{{\rm sgn \,}}
\newcommand{\rint}{\int^{\infty}_{-\infty}}
\newcommand{\dr}{\mathrm{d}}
\newcommand{\red}[1]{\textcolor{red}{#1}}
\newcommand{\skakko}[1]{\left(#1\right)}
\newcommand{\mkakko}[1]{\left\{#1\right\}}
\newcommand{\lkakko}[1]{\left[#1\right]}
\newcommand{\Z}{\field{Z}}
\newcommand{\Zo}{\field{Z}_0}
\newcommand{\abs}[1]{\lvert#1\rvert}
\newcommand{\ct}[1]{\langle#1\rangle}
\newcommand{\inp}[2]{\langle#1,#2\rangle}
\newcommand{\norm}[1]{\lVert#1 \rVert}
\newcommand{\Bnorm}[1]{\Bigl\lVert#1\Bigr  \rVert}
\newcommand{\Babs}[1]{\Bigl \lvert#1\Bigr \rvert} 
\newcommand{\ep}{\epsilon} 
\newcommand{\sumn}[1][i]{\sum_{#1 = 1}^T}
\newcommand{\tsum}[2][i]{\sum_{#1 = -#2}^{#2}}
\providecommand{\abs}[1]{\lvert#1\rvert}
\providecommand{\Babs}[1]{\Bigl \lvert#1\Bigr \rvert} 
\newcommand{\uint}{\int^{1}_{0}}
\newcommand{\freqint}{\int^{\pi}_{-\pi}}
\newcommand{\li}[1]{\mathfrak{L}(S_{#1})}
\newcommand{\cum}{{\rm cum}}
\newcommand{\xt}{\bm{X}_{t, T}}
\newcommand{\yt}{\bm{Y}_{t, T}}
\newcommand{\zt}{\bm{Z}_{t, T}}
\newcommand{\gcu}{{\rm GC}^{2 \to 1}(u)}
\newcommand{\btheta}{\bm{\bm{\theta}}}
\newcommand{\bbeta}{\bm{\eta}}
\newcommand{\bzeta}{\bm{\zeta}}
\newcommand{\bzero}{\bm{0}}
\newcommand{\bI}{\bm{I}}
\newcommand{\bd}{\bm{d}}
\newcommand{\bx}[1]{\bm{X}_{#1, T}}
\newcommand{\be}{\bm{\ep}}
\newcommand{\bp}{\bm{\phi}}
\newcommand{\act}{A_T^{\circ}}
\newcommand{\ac}{A^{\circ}}
\newcommand{\dlim}{\xrightarrow{d}}
\newcommand{\plim}{\rightarrow_{P}}
\newcommand{\ls}{\mathcal{S}}
\newcommand{\cs}{\mathcal{C}}
\newcommand{\anypath}{\,\,\textbf{-\,-\,-}\,\,}
\providecommand{\ttr}[1]{\textcolor{red}{ #1}}
\providecommand{\ttb}[1]{\textcolor{blue}{ #1}}
\providecommand{\ttg}[1]{\textcolor{green}{ #1}}
\providecommand{\tty}[1]{\textcolor{yellow}{ #1}}
\providecommand{\tto}[1]{\textcolor{orange}{ #1}}
\providecommand{\ttp}[1]{\textcolor{purple}{ #1}}
\newcommand{\sign}{\mathop{\rm sign}}
\newcommand{\conv}{\mathop{\rm conv}}
\newcommand{\argmax}{\mathop{\rm arg~max}\limits}
\newcommand{\argmin}{\mathop{\rm arg~min}\limits}
\newcommand{\argsup}{\mathop{\rm arg~sup}\limits}
\newcommand{\arginf}{\mathop{\rm arg~inf}\limits}
\newcommand{\diag}{\mathop{\rm diag}}
\newcommand{\minimize}{\mathop{\rm minimize}\limits}
\newcommand{\maximize}{\mathop{\rm maximize}\limits}
\newcommand{\tr}{\mathop{\rm tr}}
\newcommand{\Cum}{\mathop{\rm Cum}\nolimits}
\newcommand{\Var}{\mathop{\rm Var}\nolimits}
\newcommand{\Cov}{\mathop{\rm Cov}\nolimits}
\numberwithin{equation}{section}
\theoremstyle{plain}
\newtheorem{thm}{Theorem}[section]
\newtheorem{lem}[thm]{Lemma}
\newtheorem{prop}[thm]{Proposition}
\theoremstyle{definition}
\newtheorem{defi}[thm]{Definition}
\newtheorem{assumption}[thm]{Assumption}
\newtheorem{cor}[thm]{Corollary}
\newtheorem{rem}[thm]{Remark}
\newtheorem{eg}[thm]{Example}
\title{A Unified Graphical Criterion for Characterizing the Causal Interpretation of Partial Regression Coefficients in Linear Structural Equation Models}
\author{Masato Shimokawa}
\affil{Independent Researcher, Nagano, Japan}
\date{}
\maketitle
\begin{abstract}
This paper provides a graph-based characterization of partial regression coefficients in linear structural equation models.
First, we derive a generalized graphical criterion that unifies the d-separation, single-door, and back-door criteria. 
This criterion provides a generically necessary and sufficient condition under which a partial regression coefficient coincides with a linear causal effect that is not mediated by other explanatory variables. 
Second, we clarify the mechanism underlying post-treatment bias and provide a quantitative characterization of this bias.
This characterization offers a unified framework for analyzing graph structures that induce post-treatment bias, which have previously been studied on a case-by-case basis.
These results are derived from the algebraic properties of acyclic directed mixed graphs and do not rely on any specific probability distribution. 
Consequently, they apply to a broad class of linear structural equation models.
\end{abstract}

\section{Introduction}\label{sec:introduction}
Partial regression coefficients are among the most widely used quantities in empirical research, and they are often interpreted as representing the effect of a variable while holding other explanatory variables (covariates) fixed.
Nevertheless, this conventional interpretation is not generally valid because partial regression coefficients are statistical quantities rather than causal parameters.
To address this issue, the conditions under which partial regression coefficients admit a causal interpretation have been studied extensively.
For example, in linear structural equation models (SEMs) \citep{Wright1921, Haavelmo1943, Joreskog1970, Goldberger1972, Bollen1989}, it is well known that the coefficient obtained by adjusting for covariates satisfying the back-door criterion \citep{Pearl1993, Pearl1998}, or equivalently the conditional ignorability assumption \citep{Rubin1974, Angrist2009}, identifies the total effect.

The back-door criterion excludes mediators from the adjustment set and therefore does not address the causal interpretation of a partial regression coefficient when other explanatory variables include mediators.
As special cases of such situations, it is known that the d-separation and single-door criteria \citep{Spirtes1998, Pearl1998} characterize the identifiability of zero effects and direct effects, respectively.
However, these results have not been unified within a more general framework.
Our first objective is to derive a general condition under which partial regression coefficients admit a causal interpretation, even in the presence of mediator adjustment, and thereby clarify the precise scope in which the conventional interpretation is valid.

\begin{figure}[t]
\centering
\includegraphics{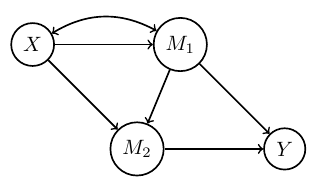}
\caption{An example where the back-door criterion is not applicable}
\label{fig:graph1}
\end{figure}

A motivating example illustrating why the conventional interpretation is important is shown in Figure~\ref{fig:graph1}.
Consider the following linear regression equation at the population level:
\begin{align*}\label{eq:intro}
Y &= \beta_{Y0\mid M_1X} + \beta_{YM_1\mid X}M_1 + \beta_{YX\mid M_1}X + \epsilon_{Y\mid M_1X},
\end{align*}
where the residual $\epsilon_{Y\mid M_1X}$ is uncorrelated with each explanatory variable.
A simple extension of the single-door criterion shows that the coefficient $\beta_{YX\mid M_1}$ coincides with the effect that is not mediated by $M_1$, that is, the effect transmitted along the specific path $X \rightarrow M_2 \rightarrow Y$.
Figure~\ref{fig:graph1} illustrates a case in which no covariate set satisfying the back-door criterion exists, rendering the total effect of $X$ on $Y$ unidentifiable, while an effect transmitted through a selected causal pathway remains identifiable. 
Thus, even when the total effect is not identifiable, useful causal information may still be obtained. 
In addition, identifying such effects can reveal which causal pathways contribute most to the overall effect.

Although Figure~\ref{fig:graph1} can be handled by a simple extension of the single-door criterion, more general settings may involve bias-inducing structures that fall outside its scope. 
For example, consider the graph obtained from Figure~\ref{fig:graph1} by removing the edges $X \leftrightarrow M_1$ and $X \rightarrow M_2$. One might expect that adjusting only for $M_2$ would identify the effect transmitted through $X \rightarrow M_1 \rightarrow Y$, because all other paths between $X$ and $Y$ are blocked. 
However, unintended post-treatment bias arises despite the absence of any unblocked path other than the target causal path.
This example demonstrates that the intuition underlying the single-door criterion does not generally extend to more complex settings.
The main difficulty in causally interpreting partial regression coefficients under adjustment for mediators lies in the possibility of unintended post-treatment bias.
Although post-treatment bias has been systematically discussed by \cite{Elwert2014} and \cite{Cinelli2024}, a detailed characterization of how such bias affects partial regression coefficients has not yet been fully developed.

The main contributions of this paper are the derivation of a general condition under which partial regression coefficients admit a causal interpretation even in the presence of mediator adjustment and the elucidation of the general mechanism underlying post-treatment bias.
First, we propose a generalized graphical criterion that unifies the d-separation, single-door, and back-door criteria. 
Based on the proposed criterion, we provide a generically necessary and sufficient condition under which a partial regression coefficient coincides with a linear causal effect that is not mediated by other explanatory variables.
This condition reveals that, contrary to the intuition underlying the single-door criterion, blocking all paths other than the target causal paths is not generally sufficient.
Second, we provide a quantitative characterization of post-treatment bias in partial regression coefficients.
This characterization reveals that unintended post-treatment bias arises from the interplay between directed paths and back-door paths.

To clarify the position of this paper, we review the studies most closely related to it.
Regression models that include mediators have long been used in path analysis and mediation analysis. 
In classical path analysis, causal effects are decomposed into direct and indirect components by means of linear SEMs, and regression coefficients are often interpreted as path coefficients representing effects conditional on other variables \citep{Wright1921, Alwin1975, Bollen1987}. 
A similar regression-based idea underlies the classical approach to mediation analysis, in which the coefficient of the treatment in a regression that includes a mediator is interpreted as a direct effect, or as the effect of the treatment while holding the mediator fixed \citep{JuddKenny1981, BaronandKenny1986}. 
\cite{Spirtes1998} and \cite{Pearl1998} established the single-door criterion, which provided a rigorous theoretical foundation for interpreting certain regression coefficients as direct effects.

At the same time, researchers have long cautioned that adjustment for post-treatment variables can be problematic.
In particular, previous work has pointed out that adjusting for mediators may not only remove mediated effects but also induce collider bias \citep{Pearl1995, Elwert2014}.
Furthermore, \citet{Cinelli2024} showed that post-treatment controls that induce unintended bias are not limited to colliders and provided a systematic overview of various forms of post-treatment bias.
By providing a unified characterization of post-treatment bias in partial regression coefficients, the results of this paper extend the findings of these previous studies.

The remainder of this paper is organized as follows.
Section~\ref{sec:preliminaries} introduces key terms as preliminaries.
Section~\ref{sec:CaR} formulates the problem addressed in this paper and presents the main theorems.
This section also discusses the mechanisms underlying post-treatment bias.
Section~\ref{sec:discussion} discusses the relationship between the proposed criterion and existing identification methods for nonparametric path-specific effects \citep{Avin2005, Shpitser2013} and linear causal parameters.
Section~\ref{sec:conc} concludes this paper.
All proofs are presented in the Appendix.

\section{Preliminaries}\label{sec:preliminaries}
\subsection{Model Setup}
Let $\bm{V} = \{X_1, \dots, X_N\}$ be a finite set of random variables with nonzero variances.
Note that $\bm{V}$ may include unobservable variables.
This paper considers the following model for the data-generating process:
\begin{equation}\label{eq:SEM}
X_i = \sum_{X_j \in \bm{PA}(X_i)} a_{ij} X_j + u_i, \quad X_i \in \bm{V},
\end{equation}
where $\bm{PA}(X_i) \subset \bm{V} \setminus \{X_i\}$ denotes the set of immediate causes of $X_i$.
Here, $a_{ij}$ is a nonzero constant referred to as the path coefficient from $X_j$ to $X_i$, and $u_1, \dots, u_N$ are error terms representing aggregations of omitted causes.
We do not assume any specific distribution for the error terms.

Model~\eqref{eq:SEM} can be rewritten as follows:
\begin{equation*}
X_i = c_i + \sum_{X_j \in \bm{PA}(X_i)} a_{ij} X_j + \tilde{u}_i, \quad X_i \in \bm{V},
\end{equation*}
where $c_i = \E[u_i]$ and $\tilde{u}_i = u_i - \E[u_i]$.
Moreover, the following vector form is useful:
\begin{equation}\label{eq:Vector SEM}
\bm{X} = \bm{A}\bm{X} + \bm{u} = \bm{c} + \bm{A}\bm{X} + \tilde{\bm{u}},
\end{equation}
where $\bm{X} = (X_1, \dots, X_N)^\top$, $\bm{u} = (u_1, \dots, u_N)^\top$, and $\tilde{\bm{u}} = (\tilde{u}_1, \dots, \tilde{u}_N)^\top$.
Here, $\bm{A}$ is an $N \times N$ matrix, and $\bm{c}$ is an $N \times 1$ vector.
We make the following assumptions on~\eqref{eq:Vector SEM}.
\begin{assumption}\label{assump:regularity DAG}
\quad
\begin{itemize}
\item[(i)] $\bm{\Sigma} \coloneq \Var(\bm{u})$ is a positive definite matrix.
\item[(ii)] $\bm{A}$ is a lower triangular matrix with all diagonal elements equal to zero.
\end{itemize}
\end{assumption}
Assumption~\ref{assump:regularity DAG} ensures that the graph associated with model~\eqref{eq:SEM} is an acyclic directed mixed graph (ADMG); see \citep{Richardson2003}.
With respect to the coefficient matrix $\bm{A}$, strictly speaking, it is sufficient that $\bm{A}$ can be transformed into the structure described in (ii) by reordering the indices of the random variables in $\bm{V}$.

\subsection{Graph Terminology and Notation}
An ADMG is a triple $(\bm{V}, \bm{E}, \bm{B})$, where $\bm{V}$ is a finite set of vertices,
$\bm{E}$ is a set of directed edges containing no directed cycles,
and $\bm{B}$ is a set of bidirected edges.
For vertices $a, b \in \bm{V}$, if they are connected by either a directed edge
(denoted by $a \rightarrow b$) or a bidirected edge (denoted by $a \leftrightarrow b$),
then $a$ and $b$ are said to be adjacent.
If there is a directed edge from $a$ to $b$, then $a$ is called a parent of $b$,
and $b$ is called a child of $a$.
A path is a sequence of distinct vertices such that consecutive vertices are adjacent.
A directed path from $a$ to $b$ (denoted by $a \rightarrow\rightarrow b$) is a path
along which all edges are directed and oriented from $a$ toward $b$.
A directed cycle is a concatenation of a directed path $a \rightarrow\rightarrow b$ and a directed edge $b \rightarrow a$.
If there is a directed path from $a$ to $b$, then $a$ is called an ancestor of $b$,
and $b$ is called a descendant of $a$.
We denote the sets of ancestors and descendants of $a$ by
$\bm{AN}(a)$ and $\bm{DE}(a)$, respectively.
A back-door path from $a$ to $b$ (denoted by $a \dashleftarrow b$)
is a path from $a$ to $b$ that begins with either a bidirected edge
or a directed edge pointing into $a$.
A v-structure is a path of the form
$a \rightarrow c \leftarrow b$,
$a \rightarrow c \leftrightarrow b$,
$a \leftrightarrow c \leftarrow b$, or
$a \leftrightarrow c \leftrightarrow b$,
where $c$ is called a collider.

\subsection{Causal Path Diagram and the Selective-Door Criterion}
We define the graph under consideration as follows.
\begin{defi}[Causal path diagram]\label{def:causal diagram}
Consider model~\eqref{eq:SEM}.
Let $\bm{V}$ denote the set of vertices, and let $\bm{E}$ denote the set of directed edges such that, for each $X_i \in \bm{V}$, the set $\bm{PA}(X_i)$ coincides with the parent set of $X_i$.
In addition, let $\bm{B}$ denote the set of bidirected edges given by
\begin{align*}
    \bm{B} \coloneq\{ \text{$X_i \leftrightarrow X_j \mid \Cov(u_i,u_j) \ne 0$ } \} .
\end{align*}
Then, the triple $(\bm{V}, \bm{E}, \bm{B})$ is called the \textit{causal path diagram} associated with model~\eqref{eq:SEM}.
\end{defi}

Next, we introduce the following concepts to visually identify the conditional relationships between variables.

\begin{defi}[Blocking]\label{def:block}
Let $\bm{Z}$ be a set such that $\bm{Z} \cap \{X_i, X_j\} = \emptyset$.
The set $\bm{Z}$ is said to \textit{block} a path $P$ between $X_i$ and $X_j$ if it satisfies one of the following two conditions; otherwise, $\bm{Z}$ is said to \textit{open} the path:
\begin{itemize}
\item[(i)] $\bm{Z}$ contains a vertex on $P$ that is not a collider in any v-structure along $P$; or
\item[(ii)] There exists a v-structure along $P$ such that neither the collider nor any of its descendants belongs to $\bm{Z}$.
\end{itemize}
\end{defi}

\cite{Pearl1993} focused on back-door paths, which can bias causal relationships, and established the back-door criterion.
This criterion for covariate selection provides a sufficient condition for identifying total causal effects.

\begin{defi}[Back-door criterion]\label{def:back-door criterion}
    A set $\bm{Z}$ satisfying $\bm{Z} \cap \{X_i,X_j\}=\emptyset$ is said to satisfy the \textit{back-door criterion} relative to an ordered pair $(X_j,X_i)$ if
\begin{enumerate}
    \item[(i)]  $\bm{Z} \cap \bm{DE}(X_j) = \emptyset$, and
    \item[(ii)] $\bm{Z}$ blocks every back-door path $X_j \dashleftarrow X_i$.
\end{enumerate}
\end{defi}

Although the back-door criterion serves as a powerful framework for covariate selection, it is not designed to capture effects transmitted through specific paths.
To address this issue, we relax the first condition of the back-door criterion.

\begin{defi}[Selective-door criterion]\label{def:selective-door criterion}
    A set $\bm{Z}$ satisfying $\bm{Z} \cap \{X_i,X_j\}=\emptyset$ is said to satisfy the \textit{selective-door criterion} relative to an ordered pair $(X_j,X_i)$ if
\begin{enumerate}
    \item[(i)] for any $X_k \in \bm{Z} \cap \bm{DE}(X_j)$ such that there exists a directed path $X_j \rightarrow\rightarrow X_k$ that is opened by $\bm{Z}\setminus\{X_k\}$,
    every back-door path $X_k \dashleftarrow X_i$ is blocked by $(\bm{Z}\cup\{X_j\})\setminus\{X_k\}$, and
    \item[(ii)] $\bm{Z}$ blocks every back-door path $X_j \dashleftarrow X_i$.
\end{enumerate}
\end{defi}
\noindent
Since the first condition is weaker than that of the back-door criterion, the selective-door criterion allows some directed paths from $X_j$ to $X_i$ to be blocked.
It should be noted, however, that blocking every path from $X_j$ to $X_i$ except for the directed paths of interest is not sufficient. 
\begin{figure}[H]
\centering
\includegraphics{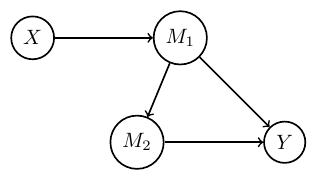}
\caption{An example where a non-collider induces unintended post-treatment bias} 
\label{fig:graph2}
\end{figure}
In Figure~\ref{fig:graph2}, the set $\{M_2\}$ leaves open only the directed path $X \rightarrow M_1 \rightarrow Y$.
However, this set does not satisfy the selective-door criterion relative to $(X,Y)$ due to the existence of the sequence of vertices $X \rightarrow M_1 \rightarrow M_2 \leftarrow M_1 \rightarrow Y$.
Since this structure is not a path and 
$M_2$ is not a collider, it can be easily overlooked when relying solely on the conventional concept of path-blocking.
In a later section, we will examine in detail how this structure induces unintended post-treatment bias.

\subsection{Controlled Total Effect}
In the literature on linear SEMs, the decomposition of causal effects has been extensively discussed and forms the foundation of mediation analysis.
In this framework, the direct effect of $X_j$ on $X_i$ is given by the path coefficient $a_{ij}$, which corresponds to the directed edge $X_j \to X_i$.
The indirect effect of $X_j$ on $X_i$ is defined as the sum of products of path coefficients along each directed path $X_j \rightarrow\rightarrow X_i$, excluding the direct edge $X_j \to X_i$.
The total effect of $X_j$ on $X_i$ is the sum of the direct and indirect effects.
See \cite{Alwin1975} and \cite{Bollen1987} for details on these causal effects.

To clarify the interpretation of partial regression coefficients, we define the following causal effect.

\begin{defi}[Controlled total effect]\label{def:controlled total effect}
Let $G = (\bm{V}, \bm{E}, \bm{B})$ be the causal path diagram associated with model~\eqref{eq:SEM}.
For a vertex set $\bm{Z} \subset \bm{V} \setminus \{X_i, X_j\}$,
consider the subgraph $G'$ obtained by removing all edges that have an arrowhead at a vertex in $\bm{Z}$.
Then, the $\bm{Z}$-controlled total effect of $X_j$ on $X_i$, denoted by $\tau_{ij\mid do(\bm{Z})}$, is defined as the total effect of $X_j$ on $X_i$ restricted to the subgraph $G'$.
\end{defi}
\noindent

The $\bm{Z}$-controlled total effect is the sum of all effects not mediated by variables in $\bm{Z}$.
The validity of this measure is justified by the following alternative expression of the data-generating process, which is derived solely through successive substitution operations in the linear structural equations.
\begin{prop}[Ancestral expansion]\label{prop:ancestral expansion}
Assume that (ii) in Assumption~\ref{assump:regularity DAG} holds. 
Let $G=(\bm{V},\bm{E}, \bm{B})$ be the causal path diagram associated with model~\eqref{eq:SEM}.
Then, for any subset $\bm{S} \subset \bm{V}\setminus \{X_i\}$, the following holds:
\begin{equation}\label{eq:ancestral expansion}
X_i = \sum_{X_j \in \bm{S}} \tau_{ij \mid do(\bm{S}\setminus \{X_j\}) } X_j + \sum_{X_k \in \bm{AN}(X_i) \setminus \bm{S}} \tau_{ik\mid do(\bm{S})}u_k + u_i.
\end{equation}
\end{prop}
If none of the omitted variables that constitute the error terms $u_k$ and $u_i$ are affected by the intervention on $X_j$, the following equality is justified:
\begin{align*}\label{eq:CTE}
    \tau_{ij\mid do(\bm{Z})}=\frac{\partial}{\partial x_j} \E[X_i\mid do(\bm{Z}=\bm{z}),\, do(X_j=x_j)],
\end{align*}
where $\bm{Z}=\bm{S}\setminus\{X_j\}$, and ``$do$'' denotes the \textit{do-operator} \citep{Pearl2009}.
Note that the notation involving the do-operator denotes probabilistic concepts of the post-intervention distribution, i.e., the distribution obtained by replacing the right-hand side of the corresponding structural equation with a constant.

\section{Causation and Regression}\label{sec:CaR}
\subsection{Problem Formulation}
Suppose we are interested in the causal effect of $X_j$ on $X_i$.
Let $\bm{S} \subset \bm{V} \setminus\{X_i\}$ denote the set of explanatory variables, and consider the following linear regression equation at the population level:
\begin{equation}\label{eq:PLRE}
X_i = \beta_{i0\mid \bm{S}} + \beta_{ij\mid \bm{Z}}X_j + \sum_{X_k\in\bm{Z}}\beta_{ik\mid \bm{S}\setminus\{X_k\}}X_k + \epsilon_{i\mid \bm{S}},
\end{equation}
where $\bm{Z} = \bm{S}\setminus\{X_j\}$ denotes a set of covariates, and the coefficients are projection coefficients onto the space spanned by explanatory variables.
See the Appendix for details.
The objective is to derive a sufficient condition on $\bm{Z}$ under which $\beta_{ij\mid \bm{Z}} = \tau_{ij\mid do(\bm{Z})}$ holds.
If the coefficient $\beta_{ij\mid \bm{Z}}$ corresponds to the causal effect, then, under standard regularity conditions for random explanatory variables, the least squares estimator provides a consistent estimate of this effect.

\subsection{Main Result}
The selective-door criterion provides a sufficient condition for this equality. 

\begin{thm}[Selective-door adjustment]\label{thm:selective-door criterion}
Assume that Assumption~\ref{assump:regularity DAG} holds. Let $G=(\bm{V},\bm{E},\bm{B})$ be the causal path diagram associated with model~\eqref{eq:SEM}. For any set $\bm{Z} \subset \bm{V}\setminus\{X_i,X_j\}$ satisfying the selective-door criterion relative to $(X_j,X_i)$, the coefficient $\beta_{ij\mid \bm{Z}}$ in \eqref{eq:PLRE} coincides with the $\bm{Z}$-controlled total effect $\tau_{ij\mid do(\bm{Z})}$.
\end{thm}

As will be shown later, the converse of Theorem~\ref{thm:selective-door criterion} also holds almost everywhere in the parameter space. That is, in a generic sense, the selective-door criterion constitutes a necessary and sufficient condition for providing a linear causal interpretation of partial regression coefficients.

As corollaries, we obtain the following propositions.

\begin{cor}[D-separation adjustment]\label{col:DSA}
Assume that Assumption~\ref{assump:regularity DAG} holds. Let $G=(\bm{V},\bm{E},\bm{B})$ be the causal path diagram associated with model~\eqref{eq:SEM}. For any set $\bm{Z} \subset \bm{V}\setminus \{X_i,X_j\}$ that blocks every path between $X_j$ and $X_i$, the coefficient $\beta_{ij\mid \bm{Z}}$ in \eqref{eq:PLRE} is equal to $0$.
\end{cor}

\begin{cor}[Single-door adjustment]\label{col:SDA}
Assume that Assumption~\ref{assump:regularity DAG} holds. Let $G=(\bm{V},\bm{E},\bm{B})$ be the causal path diagram associated with model~\eqref{eq:SEM}. For any set $\bm{Z} \subset \bm{V}\setminus (\{X_i,X_j\}\cup\bm{DE}(X_i))$ that blocks every path between $X_j$ and $X_i$ except for the directed edge $X_j \rightarrow X_i$, the coefficient $\beta_{ij\mid \bm{Z}}$ in \eqref{eq:PLRE} coincides with the direct effect of $X_j$ on $X_i$.
\end{cor} 

\begin{cor}[Back-door adjustment]\label{col:BDA}
Assume that Assumption~\ref{assump:regularity DAG} holds. Let $G=(\bm{V},\bm{E},\bm{B})$ be the causal path diagram associated with model~\eqref{eq:SEM}. For any set $\bm{Z} \subset \bm{V}\setminus (\{X_i,X_j\}\cup\bm{DE}(X_j))$ satisfying the back-door criterion relative to $(X_j,X_i)$, the coefficient $\beta_{ij\mid \bm{Z}}$ in \eqref{eq:PLRE} coincides with the total effect of $X_j$ on $X_i$.
\end{cor}

These corollaries were originally established in \cite{Spirtes1998} and \cite{Pearl1998}.
Note that the single-door criterion (adjustment) excludes descendants of the outcome variable.
The reason for this exclusion is explained by the mechanism of post-treatment bias discussed in the following subsection.

In many empirical analyses, it is not uncommon for all explanatory variables to be regarded as causes of interest. 
Thus, we also provide sufficient conditions under which all regression coefficients can be interpreted as causal effects.

\begin{prop}\label{prop:no confundings}
The following statements are equivalent:
\begin{itemize}
    \item[(i)] For any $X_j \in \bm{S}$, the set $\bm{S}\setminus\{X_j\}$ satisfies the selective-door criterion relative to $(X_j,X_i)$.
    \item[(ii)] For any $X_j \in \bm{S}$, the set $\bm{S}\setminus\{X_j\}$ blocks every back-door path $X_j \dashleftarrow X_i$.
    \item[(iii)] For any $X_j \in \bm{S}$, the set $\bm{S}\setminus\{X_j\}$ blocks every back-door path $X_j \dashleftarrow X_i$ that contains no v-structure.
\end{itemize}
\end{prop}
\noindent
More precisely, Proposition \ref{prop:no confundings} establishes a sufficient condition for the vector of partial regression coefficients to coincide with the total joint effects \citep{Nandy2017}.
Condition $(iii)$ is particularly easy to check because $\bm{S}\setminus\{X_j\}$ blocks a back-door path from $X_j$ to $X_i$ that contains no v-structure if and only if $X_i \notin \bm{S}$ and at least one vertex $X_k \in \bm{S}\setminus\{X_j\}$ lies on that back-door path.

\subsection{Post-Treatment Bias}
Biases that arise from conditioning on post-treatment variables (i.e., descendants in the graph) are  collectively referred to as post-treatment bias.
Examining the mechanism of this bias helps clarify why the first condition of the selective-door criterion is required.

In addition to the set $\bm{S}=\{X_j\} \cup \bm{Z}$ in \eqref{eq:PLRE}, the following disjoint subsets are introduced:
\begin{align*}
    \bm{S}_1 &= \{\, X_k \in \bm{S} \mid \text{$\bm{S} \setminus \{X_k\}$ opens a back-door path $X_k \dashleftarrow X_i$}. \,\}, \\
    \bm{S}_2 &= \{\, X_k \in \bm{S} \mid \text{$\bm{S} \setminus \{X_k\}$ blocks every back-door path $X_k \dashleftarrow X_i$}. \,\}.
\end{align*}

The bias $\gamma_{ik\mid \bm{S}\setminus\{X_k\}}$ corresponding to each coefficient $\beta_{ik\mid \bm{S}\setminus\{X_k\}}$ in \eqref{eq:PLRE} is defined as
\[
\gamma_{ik\mid \bm{S}\setminus\{X_k\}} \coloneq \beta_{ik\mid \bm{S}\setminus\{X_k\}} - \tau_{ik\mid do(\bm{S}\setminus\{X_k\})}, \quad X_k \in \bm{S}.
\]
Although bias is often defined as the difference from the total effect, we adopt the above definition for clarity in the quantitative characterization of post-treatment bias.

The following theorem characterizes post-treatment bias after confounding bias has been controlled by blocking back-door paths.

\begin{thm}\label{thm:PTB}
Assume that Assumption~\ref{assump:regularity DAG} holds. Let $G=(\bm{V},\bm{E},\bm{B})$ be the causal path diagram associated with model~\eqref{eq:SEM}. For any subset $\bm{Z} \subset \bm{V}\setminus\{X_i,X_j\}$ that blocks every back-door path $X_j \dashleftarrow X_i$, the bias of the coefficient $\beta_{ij\mid \bm{Z}}$ in \eqref{eq:PLRE} is given by
\[
\gamma_{ij\mid \bm{Z}} = -\sum_{X_p\in\bm{S}_1}\gamma_{ip\mid \bm{S}\setminus\{X_p\}}\tau_{pj\mid do(\bm{S}_2\setminus\{X_j\})}.
\]
\end{thm}

The first condition of the selective-door criterion implies that $\tau_{pj\mid do(\bm{S}_2\setminus\{X_j\})}=0$ for the relevant vertices $X_p\in\bm{S}_1$ and can be understood as a condition to prevent unintended post-treatment bias. 
If the total effect is the causal parameter of interest, the following decomposition is useful.
\[
\beta_{ij\mid \bm{Z}} = \tau_{ij} - ( \tau_{ij} - \tau_{ij|do(\bm{Z})} ) + \gamma_{ij\mid \bm{Z}},
\]
where $\tau_{ij} - \tau_{ij|do(\bm{Z})}$ is the mediated effect through $\bm{Z}$.
The negative second term represents the component of post-treatment bias that arises from blocking mediated pathways, whereas the third term represents unintended post-treatment bias arising from mechanisms other than mediated-path blocking.

\begin{figure}[t]
    \centering
    \includegraphics[scale=1.0]{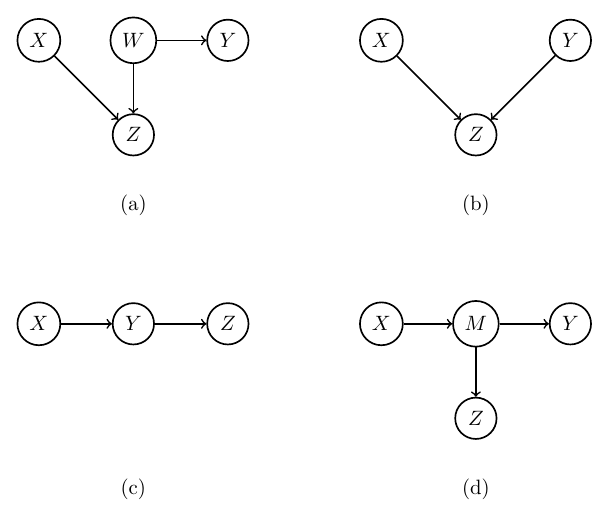}
    \caption{Graph structures inducing unintended post-treatment bias}
    \label{fig:fig3}
\end{figure}

Figure~\ref{fig:fig3} illustrates four well-known structures that induce unintended post-treatment bias.
In cases (a) and (b), as considered under path-blocking, including only the collider $Z$ as a covariate leads to bias.
Case (c) corresponds to a virtual-collider-type (case-control-type) structure discussed by \cite{Pearl2009} and \cite{Cinelli2024}, whereas case (d) arises from conditioning on a descendant of a mediator.
Although these structures are often discussed separately, Theorem~\ref{thm:PTB} shows that they can all be described in a unified way as involving both open directed paths and open back-door paths.
In particular, while case (d) is treated as an overcontrol-type bad control by \cite{Cinelli2024}, classifying this structure as overcontrol bias may be somewhat misleading, at least within the class of linear models.
Theorem~\ref{thm:PTB} suggests that this bias is not best understood as a consequence of indirectly controlling for the mediator; rather, it arises from the same general mechanism as cases (a), (b), and (c), namely the interplay between open directed paths and open back-door paths.

The reason the single-door criterion excludes descendants of the outcome variable from the set of covariates is that it accounts for the possibility of structures such as (c).
Note that this criterion assumes the blocking of all directed paths except for the directed edge; hence, structures such as (d) need not be considered.
Furthermore, under the d-separation criterion, all directed paths are blocked; hence, including descendants of the outcome variable does not induce a nonzero partial regression coefficient for a zero target effect.

\subsection{Necessity for Almost All Parameter Points}
Here, we discuss the precise scope of identifiability of the selective-door criterion.
For the coefficient matrix $\bm{A}$ and covariance matrix $\bm{\Sigma}$, let $(\bm{A},\bm{\Sigma})$ denote a parameter of the model~\eqref{eq:SEM}.
Since covariance information is sufficient to characterize the partial regression coefficients of interest, we omit the intercept term $\bm{c}$ from the model parameter.
We consider the following parameter space:
\begin{align}\label{eq:parameter space}
    \bm{\Theta} \coloneq \{ (\bm{A},\bm{\Sigma}) \mid \text{ $\bm{A}$ and $\bm{\Sigma}$ satisfy Assumption  2.1} \}.
\end{align}
Almost everywhere in this parameter space, satisfying the selective-door criterion is not only sufficient but also necessary for partial regression coefficients to admit a linear causal interpretation.
\begin{thm}\label{thm:Necessity almost everywhere}
Let $G=(\bm{V},\bm{E},\bm{B})$ be the causal path diagram associated with model~\eqref{eq:SEM}.
Fix a set $\bm{Z} \subset \bm{V} \setminus \{X_i, X_j\} $ and consider the coefficient $\beta_{ij \mid \bm{Z}}$ in \eqref{eq:PLRE}.
Then, except on a set of Lebesgue measure zero in the parameter space $\bm{\Theta}$, the following statements are equivalent:
\begin{enumerate}
    \item[(i)] The set $\bm{Z}$ satisfies the selective-door criterion relative to $(X_j, X_i)$.
    \item[(ii)] The coefficient $\beta_{ij \mid \bm{Z}}$ coincides with the $\bm{Z}$-controlled total effect $\tau_{ij \mid do(\bm{Z})}$.
\end{enumerate}
\end{thm}

\section{Discussion}\label{sec:discussion}

In this section, we compare the proposed criterion with existing identification results for nonparametric path-specific effects and linear causal parameters.

\subsection{Relationship to Nonparametric Path-Specific Effect Identification}

\cite{Avin2005} introduced the recanting witness criterion, and \cite{Shpitser2013} generalized it and proposed the recanting district criterion, a graphical condition for the identifiability of path-specific effects in nonparametric causal models. 
These criteria characterize whether path-specific effects are identifiable as functionals of ordinary interventional distributions in nonparametric causal models. 
The parameter considered in this paper concerns the causal interpretation of partial regression coefficients in linear structural equation models, and therefore differs from the nonparametric path-specific effects studied in their work. 
Nevertheless, the two frameworks are closely related in that both aim to separate effects transmitted through selected causal pathways from those transmitted through other pathways.

An important point in understanding this relationship is that the recanting district criterion and the selective-door criterion are not formulated under the same assumptions. 
The recanting district criterion addresses whether a nonparametric path-specific effect is identifiable, or more precisely, whether it can be represented in terms of interventional distributions, assuming that any vertex in the ADMG can be observed if necessary. 
In contrast, the selective-door criterion focuses on the value of a partial regression coefficient and distinguishes between vertices that must be observed and those that need not be observed in order to identify a controlled total effect. 
Thus, the latter should not be viewed as merely addressing a complete special case of the former.

Although it is not necessarily appropriate to seek a simple inclusion relationship between the two criteria, the following proposition suggests that the graphical structures captured by them share a common underlying feature.

\begin{prop}\label{prop:SDCvsRDC}
Let $\pi$ denote the collection of directed paths from $X$ to $Y$ that do not pass through any vertex in $\bm{Z}$.
If $\bm{Z}$ satisfies the first condition of the selective-door criterion relative to $(X,Y)$,
then no recanting district exists for the $\pi$-specific effect.
\end{prop}

The converse of this proposition does not hold in general. 
This is largely due to the difference in the problem settings. 
For example, Figure~\ref{fig:fig3} (c) provides a counterexample.

\subsection{Relationship to Identification Methods for Linear Causal Parameters}

A large literature has studied the identification of causal parameters in linear SEMs. The single-door and back-door criteria \citep{Spirtes1998, Pearl1998} identify direct and total effects, respectively. Extensions of instrumental-variable methods were developed by \cite{BritoandPearl2002a}, \cite{BritoandPearl2002b}, and \cite{Chan2010}, while alternative approaches were proposed by \cite{KurokiandPearl2014} and \cite{Stanghellini2015}. Furthermore, \cite{Foygel2012} introduced the half-trek criterion for the generic identification of model parameters, which was subsequently extended by \cite{Chen2016} and \cite{Barber2022}.
These studies primarily focus on the identification of direct effects, total effects, or the full set of model parameters. In contrast, the causal parameter considered in this paper is a linear path-specific effect that explicitly depends on the covariate set. Consequently, the selective-door criterion is complementary to existing identification criteria rather than being subsumed by them.

\cite{Nandy2017} studied the identification and estimation of total joint effects in linear SEMs. They pointed out that, unlike total effects of single interventions, total joint effects cannot in general be obtained by straightforward covariate adjustment from a single multiple regression, and proposed methods based on recursive regressions for causal effects and modified Cholesky decompositions.
In contrast, Proposition~\ref{prop:no confundings} addresses when the vector of partial regression coefficients itself coincides with the total joint effect.

\section{Conclusion}\label{sec:conc}
We have developed a graphical framework for characterizing precisely when partial regression coefficients admit a causal interpretation in linear structural equation models. The proposed selective-door criterion unifies the d-separation, single-door, and back-door criteria and provides a generically necessary and sufficient condition under which a partial regression coefficient coincides with a controlled total effect, that is, the effect not mediated by the other explanatory variables included in the regression model.

The results also clarify why post-treatment variable adjustment can generate bias even when all paths other than the target causal paths appear to be blocked. In particular, the analysis shows that unintended post-treatment bias arises from the interplay between directed paths and back-door paths. This provides a unified explanation of several bias-generating structures that have previously been discussed separately, including collider-type, virtual-collider-type, and descendant-of-mediator structures.
However, because the analysis in this paper is restricted to linear models, whether similar results can be established for nonparametric causal models remains an open question for future research.

Finally, the proposed results are derived from algebraic properties of acyclic directed mixed graphs and do not rely on any specific distributional assumptions. 
Since bidirected edges are used here to represent arbitrary dependence structures among error terms rather than only latent confounders with linear causal effects, the framework applies to a broad class of linear graphical models. 
These findings help delineate the precise conditions under which linear regression coefficients can be given a rigorous causal interpretation.

\section*{Acknowledgments}
The author is grateful to Professor Kou Fujimori for helpful comments and encouragement during the preparation of this manuscript.

\section*{Statements and Declarations}
The author has no relevant financial or non-financial interests to disclose.

\bibliographystyle{econ}
\bibliography{bibliography}

\section{Appendix}\label{sec:appendix}

\subsection{Partial Regression Coefficient at the Population Level}

\begin{defi}[Partial regression coefficient]\label{def:PPRC}
Let $Y$ be a random variable and $\bm{X}$ a random vector.
Suppose that the inverse matrix $\E[\bm{X}\bm{X}^\top]^{-1}$ and the vector $\E[\bm{X}Y]$ exist.
Then,
\[
\bm{\beta}_{Y\bm{X}}\coloneq\E[\bm{XX}^\top]^{-1}\E[\bm{X}Y]
\]
is called the \textit{partial regression coefficient} in the regression of $Y$ on $\bm{X}$.
\end{defi}

\begin{defi}[Linear regression equation]\label{def:PLRE}
Let $\bm{\beta}_{Y\bm{X}}$ denote the partial regression coefficient in the regression of $Y$ on $\bm{X}$.
The equation
\[
Y = \bm{X}^\top \bm{\beta}_{Y\bm{X}} + \epsilon_{Y\mid \bm{X}},
\]
where $\epsilon_{Y\mid \bm{X}} \coloneq Y - \bm{X}^\top \bm{\beta}_{Y\bm{X}}$,
is called the \textit{linear regression equation} of $Y$ on $\bm{X}$.
\end{defi}

\begin{prop}\label{prop:orthgonality}
For the residual $\epsilon_{Y\mid \bm{X}}$, the following properties hold:
\begin{enumerate}
\item[(i)] $\E[\bm{X}\epsilon_{Y\mid \bm{X}}] = \bm{0}$.
\item[(ii)] If $\bm{X}$ includes a nonzero constant c, then $\E[\epsilon_{Y\mid \bm{X}}] = 0$ and $\Cov(\bm{X}, \epsilon_{Y\mid \bm{X}}) = \bm{0}$.
\end{enumerate}
\end{prop}
\begin{proof}
$(i)$ The result follows directly from Definition~\ref{def:PPRC}.

\noindent
$(ii)$ Since $\bm{X}$ contains a nonzero constant $c$, the expectation $\E[\bm{X}\epsilon_{Y\mid \bm{X}}]$ includes
\[
\E[c\epsilon_{Y\mid \bm{X}}] = c\E[\epsilon_{Y\mid \bm{X}}].
\]
Hence, by the first statement $(i)$, we have $\E[\epsilon_{Y\mid \bm{X}}]=0$. Consequently, $\Cov(\bm{X},\epsilon_{Y\mid \bm{X}})=\bm{0}$ holds trivially.
\end{proof}

\begin{prop}\label{prop:uniqueness}
Let $Y$ be a random variable and $\bm{X}$ a random vector, and suppose that the inverse matrix $\E[\bm{X}\bm{X}^\top]^{-1}$ and the vector $\E[\bm{X}Y]$ exist.
Assume that the following relationship holds:
\[
Y = \bm{X}^\top\bm{b}+e,
\]
where $\bm{b}$ is a vector and $e$ is a random variable satisfying $\E[\bm{X}e] = \bm{0}$.
Then, it follows that $\bm{b} = \E[\bm{X}\bm{X}^\top]^{-1}\E[\bm{X}Y]$ and $e=\epsilon_{Y\mid \bm{X}}$.  
\end{prop}
\begin{proof}
By multiplying both sides of the equation by $\bm{X}$ from the left and taking expectations, we have
\[
\E[\bm{X}Y] = \E[\bm{XX}^\top]\bm{b}.
\]
Due to the existence of $\E[\bm{XX}^\top]^{-1}$, the conclusion is obtained.
\end{proof}

\begin{prop}\label{prop:the existance of PLRE}
Suppose that Assumption~\ref{assump:regularity DAG} holds.
Then, for the random vector $\bm{X}_{sub+} \coloneq (1, \bm{X}^\top_{sub})^\top$, where $\bm{X}_{sub}$ is any random subvector of $\bm{X}$ generated by equation~\eqref{eq:Vector SEM}, the inverse matrix $\E[\bm{X}_{sub+}\bm{X}_{sub+}^\top]^{-1}$ exists.
\end{prop}

\begin{proof}
Let ${\bm{X}}_{+}$ denote $(1,\bm{X}^\top)^\top$. It suffices to prove the existence of $\E[\bm{X}_{+}\bm{X}_{+}^\top]^{-1}$.
Since $\bm{I}_N-\bm{A}$, where $\bm{I}_N$ is the $N \times N$ identity matrix, is a lower unitriangular matrix and nonsingular, \eqref{eq:Vector SEM} has another form
\[
\bm{X} = \bm{B}\bm{u}=\bm{Bc}+\bm{B}\tilde{\bm{u}}, \quad \bm{B} \coloneq (\bm{I}_N-\bm{A})^{-1}.
\]
Thus, $\Var(\bm{X})=\bm{B}\bm{\Sigma}\bm{B}^\top$ is positive definite, and the expectation
\[
\E[\bm{XX}^\top] = \Var(\bm{X}) + \bm{Bcc}^\top\bm{B}^\top
\]
is also positive definite. Note that $\E[\bm{X}_{+}\bm{X}_{+}^\top]$ has the following structure:
\[
\E[\bm{X}_{+}\bm{X}_{+}^\top] =
\begin{pmatrix}
1 & \bm{c}^\top\bm{B}^\top \\
\bm{Bc} & \E[\bm{XX}^\top]
\end{pmatrix}.
\]
The conclusion follows from the properties of the determinant.
\[
\det(\E[\bm{X}_{+}\bm{X}_{+}^\top]) = 1\cdot\det(\E[\bm{XX}^\top]-\bm{Bc}
\bm{c}^\top\bm{B}^\top) = \det(\Var(\bm{X})) > 0.
\]
\end{proof}

\subsection{Proofs}
\begin{proof}[Proof of Proposition~\ref{prop:ancestral expansion}]
Equation~\eqref{eq:ancestral expansion} is obtained by iteratively substituting structural equations, where each equation has a left-hand side not belonging to $\bm{S}$, until all vertices $X_l \notin \bm{S}$ are eliminated from the right-hand side. This substitution procedure is guaranteed to terminate, since $G$ has no cycles and $\bm{V}$ is a finite set.
\end{proof}

To prove the theorems presented in this paper, we introduce some random variables and establish several lemmas.
Fix an arbitrary vertex $X_i \in \bm{V}$ and a set $\bm{S} \subset \bm{V} \setminus \{X_i\}$.
For $X_i$, each $X_p \in \bm{S}_1$, and each $X_q \in \bm{S}_2$, define the random variables
$\tilde{X}_i$, $\tilde{X}_p$, and $\tilde{X}_q$ as follows:
\begin{align*}
    \tilde{X}_i &\coloneq X_i -
    \sum_{X_k \in \bm{S}}
    \tau_{ik\mid do(\bm{S}\setminus\{X_k\})} X_k,
    \\
    \tilde{X}_p &\coloneq X_p -
    \sum_{X_q \in \bm{S}_2}
    \tau_{pq\mid do(\bm{S}_2\setminus\{X_q\})} X_q,
    \\
    \tilde{X}_q &\coloneq X_q -
    \sum_{X_p \in \bm{S}_1}
    \tau_{qp\mid do(\bm{S}_1\setminus\{X_p\})} X_p,
\end{align*}
where the disjoint subsets $\bm{S}_1, \bm{S}_2 \subset \bm{S}$ are defined as
\begin{align*}
      \bm{S}_1 &\coloneq \{\, X_k \in \bm{S}
      \mid \textit{$\bm{S} \setminus \{X_k\}$ opens a back-door path
      $X_k \dashleftarrow X_i$} \,\}, \\
      \bm{S}_2 &\coloneq \{\, X_k \in \bm{S}
      \mid \textit{$\bm{S} \setminus \{X_k\}$ blocks every back-door path
      $X_k \dashleftarrow X_i$} \,\}.
\end{align*}

\begin{lem}\label{lem:lemma1}
Assume that Assumption~\ref{assump:regularity DAG} holds. Let $G=(\bm{V}, \bm{E}, \bm{B})$ be the causal path diagram associated with model~\eqref{eq:SEM}. 
Then, $\tilde{X}_q$ is uncorrelated with $\tilde{X}_i$.
\end{lem}
\begin{proof}
Note that the following equations hold from Proposition~\ref{prop:ancestral expansion}.
\begin{align*}
    \tilde{X}_i &=\sum_{X_l \in \bm{AN}(X_i)\setminus\bm{S}} \tau_{il\mid do(\bm{S})}u_l+u_i, \\
     \tilde{X}_q &=\sum_{X_l \in \bm{AN}(X_q)\setminus\bm{S}_1} \tau_{ql\mid do(\bm{S}_1)}u_l+u_q.
\end{align*}    
If $\Cov(\tilde{X}_q,\tilde{X}_i)\ne0$ holds, there exists at least one of the following paths:
\begin{itemize}
    \item[(a)] $X_q \leftrightarrow X_i$;
    \item[(b)] $X_q \leftarrow\leftarrow X_{l} \rightarrow\rightarrow X_i$, where $X_q \leftarrow\leftarrow X_{l}$ and $X_{l} \rightarrow\rightarrow X_i$ are opened by $\bm{S}_1$ and $\bm{S}$, respectively;
    \item[(c)] $X_q \leftarrow\leftarrow X_{l_1} \leftrightarrow X_{l_2} \rightarrow\rightarrow X_i$, where $X_q \leftarrow\leftarrow X_{l_1}$ and $X_{l_2} \rightarrow\rightarrow X_i$ are opened by $\bm{S}_1$ and $\bm{S}$, respectively;
    \item[(d)] $X_q \leftarrow\leftarrow X_i$ that is opened by $\bm{S}_1$;
\end{itemize}
Path (a) contradicts the condition that $X_q \in \bm{S}_2$.
For paths (b), (c), and (d), if 
$\bm{S}_2 \setminus \{X_q\}$ were to block the subpaths that are not blocked by 
$\bm{S}_1$, then there would exist a back-door path from some $X'_q \in \bm{S}_2$ to 
$X_i$ that is opened by $\bm{S}\setminus \{X'_q\}$.
Hence, these paths are opened by $\bm{S} \setminus \{X_q\}$, which contradicts $X_q \in \bm{S}_2$.
\end{proof}

\begin{lem}\label{lem:lemma2}
Let $(\bm{V}, \bm{E}, \bm{B})$ be an ADMG.
For any vertices $X_p \in \bm{S}_1$, $X_q \in \bm{S}_2$, and
$X_l, X_{l_1}, X_{l_2} \in \bm{V} \setminus \bm{S}$,
the following paths do not exist:
\begin{enumerate}
    \item[(a)]$X_q \leftarrow\leftarrow X_l \rightarrow\rightarrow X_p$,
    where $X_q \leftarrow\leftarrow X_l$ and
    $X_l \rightarrow\rightarrow X_p$ are opened by
    $\bm{S}_1$ and $\bm{S}_2$, respectively;
    \item[(b)]$X_q \leftarrow\leftarrow X_{l_1} \leftrightarrow X_{l_2} \rightarrow\rightarrow X_p$,
    where $X_q \leftarrow\leftarrow X_{l_1}$ and
    $X_{l_2} \rightarrow\rightarrow X_p$ are opened by
    $\bm{S}_1$ and $\bm{S}_2$, respectively.
\end{enumerate}

\end{lem}

\begin{proof}
Since the nonexistence of path (b) can be proved in exactly the same way, we prove only the nonexistence of path (a).
If (a) exists, there exists the following sequence of vertices:
\[
X'_q \leftarrow\leftarrow X_l \rightarrow\rightarrow X'_p \dashleftarrow X_i,
\]
where $X'_q \in \bm{S}_2$ and $X'_p \in \bm{S}_1$ are the vertices closest to $X_l$ along the subpaths $X_q \leftarrow\leftarrow X_l$ and $X_l \rightarrow\rightarrow X_p$, respectively.
Note that this sequence of vertices satisfies the following conditions:
\begin{itemize}
    \item[(i)] $X'_q \leftarrow\leftarrow X_l$ is opened by $\bm{S}\setminus\{X'_q\}$;
    \item[(ii)] $X_l \rightarrow\rightarrow X'_p$ is opened by $\bm{S}\setminus\{X'_p\}$;
    \item[(iii)] $X'_p \dashleftarrow X_i$ is opened by $\bm{S}\setminus\{X'_p\}$; and
    \item[(iv)] $X'_p \dashleftarrow X_i$ may overlap with the path $X'_q \leftarrow\leftarrow X_l \rightarrow\rightarrow X'_p$;
\end{itemize}
This structure implies that there exists a back-door path $X'_q \dashleftarrow X_i$ opened by
$\bm{S} \setminus \{X'_q\}$, regardless of whether the paths
$X'_p \dashleftarrow X_i$ and $X'_q \leftarrow\leftarrow X_l \rightarrow\rightarrow X'_p$ overlap.
This leads to a contradiction.
\end{proof}

\begin{lem}\label{lem:lemma3}
Assume that Assumption~\ref{assump:regularity DAG} holds. Let $G=(\bm{V}, \bm{E}, \bm{B})$ be the causal path diagram associated with model~\eqref{eq:SEM}.
Then, $\tilde{X}_q$ is uncorrelated with $\tilde{X}_p$.
\end{lem}
\begin{proof}
If $\Cov(\tilde{X}_q, \tilde{X}_p)\ne0$ holds, there exists at least one of the following paths:
\begin{itemize}
    \item[(a)] $X_q \leftrightarrow X_p$;
    \item[(b)] $X_q \leftarrow\leftarrow X_{l} \rightarrow\rightarrow X_p$, where $X_q \leftarrow\leftarrow X_{l}$ and $X_{l} \rightarrow\rightarrow X_p$ are opened by $\bm{S}_1$ and $\bm{S}_2$, respectively;
    \item[(c)] $X_q \leftarrow\leftarrow X_{l_1} \leftrightarrow X_{l_2} \rightarrow\rightarrow X_p$, where $X_q \leftarrow\leftarrow X_{l_1}$ and $X_{l_2} \rightarrow\rightarrow X_p$ are opened by $\bm{S}_1$ and $\bm{S}_2$, respectively;
\end{itemize}
Paths (b) and (c) do not exist by Lemma~\ref{lem:lemma2}.
Thus, it suffices to prove that path (a) does not exist.
If (a) exists, there exists the following sequence of vertices:
\begin{align*}
    X_q \leftrightarrow X_p \dashleftarrow X_i,
\end{align*}
where
\begin{itemize}
    \item[(i)] $X_p \dashleftarrow X_i$ is opened by $\bm{S}\setminus\{X_p\}$; and
    \item[(ii)] $X_p \dashleftarrow X_i$ may overlap with the path $X_q \leftrightarrow X_p$;
\end{itemize}
This structure implies that there exists a back-door path $X_q \dashleftarrow X_i$ opened by $\bm{S} \setminus \{X_q\}$, leading to a contradiction. 
\end{proof}

\begin{proof}[Proof of Theorem~\ref{thm:selective-door criterion}]
By Proposition~\ref{prop:the existance of PLRE}, the following regression equation exists:
\[
\tilde{X}_i = \tilde{\gamma}_{i0\mid \bm{S}_1} + \sum_{X_p \in \bm{S}_1}\tilde{\gamma}_{ip\mid \bm{S}_1\setminus\{X_p\}}\tilde{X}_p + \tilde{\epsilon}_{i\mid \bm{S}_1}.
\]
As a first step, we will show that $\tilde{\epsilon}_{i\mid \bm{S}_1}=\epsilon_{i\mid \bm{S}}$ by mathematical induction.
Let the elements of 
$\bm{S}_2$ be ordered by increasing indices as $X_{q^1},...,X_{q^M}$.
Then, for $X_{q^1}$, we have
\begin{align*}
\Cov(X_{q^1},\tilde{\epsilon}_{i\mid \bm{S}_1})&=\Cov \left( \tilde{X}_{q^{1}}, \tilde{\epsilon}_{i\mid \bm{S}_1} \right)\\
&= \Cov(\tilde{X}_{q^{1}}, \tilde{X}_i) - \sum_{X_p \in \bm{S}_1}\tilde{\gamma}_{ip\mid \bm{S}_1\setminus\{X_p\}}\Cov\left(\tilde{X}_{q^{1}}, \tilde{X}_p \right) 
\end{align*}
Note that for any $X_p \in \bm{S}_1 \cap \bm{AN}(X_{q^1})$, the equality $X_p=\tilde{X}_p$ holds due to the triangularity of the coefficient matrix $\bm{A}$ and the minimality of $q^1$. 
By Lemmas~\ref{lem:lemma1} and~\ref{lem:lemma3}, the right-hand side vanishes; hence, for $q=q^{1}$, $\Cov(X_{q},\tilde{\epsilon}_{i\mid \bm{S}_1})=0$.

Next, for $q^1\leq q<q^{m}$, assume that $\Cov(X_{q},\tilde{\epsilon}_{i\mid \bm{S}_1})=0$ holds, and we prove the uncorrelatedness between $X_{q^{m}}$ and $\tilde{\epsilon}_{i\mid \bm{S}_1}$.
\begin{align*}
    \Cov(X_{q^m},\tilde{\epsilon}_{i\mid \bm{S}_1}) &= \Cov(\tilde
    {X}_{q^m},\tilde{\epsilon}_{i\mid \bm{S}_1}) + \sum_{X_p\in \bm{S}_1}\tau_{q^mp\mid do( \bm{S}_1 \setminus\{X_p\} )}\Cov(X_p,\tilde{\epsilon}_{i\mid \bm{S}_1}) \\
    &= \Cov(\tilde
    {X}_{q^m},\tilde{\epsilon}_{i\mid \bm{S}_1}) \\
    & \, + \sum_{X_p\in \bm{S}_1}\tau_{q^mp\mid do( \bm{S}_1 \setminus\{X_p\} )}\Cov(\tilde{X}_p,\tilde{\epsilon}_{i\mid \bm{S}_1}) \\
    & \, + \sum_{X_p\in \bm{S}_1}\sum_{X_q\in \bm{S}_2}
    \tau_{q^mp\mid do( \bm{S}_1 \setminus\{X_p\} )}\tau_{pq\mid do( \bm{S}_2\setminus\{X_q\} )}\Cov(X_q,\tilde{\epsilon}_{i\mid \bm{S}_1}) \\
    &= \Cov(\tilde
    {X}_{q^m},\tilde{\epsilon}_{i\mid \bm{S}_1}) \\
    & \, + \sum_{X_p\in \bm{S}_1}\sum_{X_q\in \bm{S}_2}
    \tau_{q^mp\mid do( \bm{S}_1 \setminus\{X_p\} )}\tau_{pq\mid do( \bm{S}_2\setminus\{X_q\} )}\Cov(X_q,\tilde{\epsilon}_{i\mid \bm{S}_1})
\end{align*}
Since each $X_q$ having a nonzero coefficient on the right-hand side is an ancestor of $X_{q^m}$, the covariance vanishes by assumption. 
Hence, by Lemmas~\ref{lem:lemma1} and~\ref{lem:lemma3}, the following uncorrelatedness holds.
\begin{align*}
   \Cov(X_{q^m},\tilde{\epsilon}_{i\mid \bm{S}_1}) &= \Cov(\tilde
    {X}_{q^m},\tilde{\epsilon}_{i\mid \bm{S}_1}) \\
    &= \Cov(\tilde{X}_{q^{m}}, \tilde{X}_i) - \sum_{X_p \in \bm{S}_1}\tilde{\gamma}_{ip\mid \bm{S}_1\setminus\{X_p\}}\Cov\left(\tilde{X}_{q^{m}}, \tilde{X}_p \right) \\
    &= 0
\end{align*}
Thus, $\tilde{\epsilon}_{i\mid \bm{S}_1}$ is uncorrelated with all vertices in $\bm{S}_2$, which implies that
$\tilde{\epsilon}_{i\mid \bm{S}_1}$ is also uncorrelated with all vertices in $\bm{S}_1$.
Therefore, from Proposition~\ref{prop:uniqueness}, it follows that
\begin{align*}
    \epsilon_{i\mid \bm{S}} &= \tilde{X}_i - \beta_{i0\mid \bm{S}} -\sum_{X_k \in \bm{S}}\gamma_{ik\mid \bm{S}\setminus \{X_k\}}X_k \\
    &=\tilde{X}_i - \tilde{\gamma}_{i0\mid \bm{S}_1} - \sum_{X_p \in \bm{S}_1} \tilde{\gamma}_{ip\mid \bm{S}_1\setminus\{X_p\}} \left( X_p-\sum_{X_q \in \bm{S}_2} \tau_{pq\mid do( \bm{S}_2 \setminus \{X_q\})}X_q \right)  \\
    &=\tilde{X}_i - \tilde{\gamma}_{i0\mid \bm{S}_1} - \sum_{X_p \in \bm{S}_1}\tilde{\gamma}_{ip\mid \bm{S}_1\setminus \{X_p\}}\tilde{X}_p \\
    &=\tilde{\epsilon}_{i\mid \bm{S}_1}.
\end{align*}

By Proposition~\ref{prop:uniqueness} again, we obtain the following equations.
\begin{align*}
\beta_{i0\mid \bm{S}}&= \tilde{\gamma}_{i0\mid \bm{S}_1}, \\
\gamma_{ip\mid \bm{S}\setminus\{X_p\}} &= \tilde{\gamma}_{ip\mid \bm{S}_1\setminus\{X_p\}}, \quad X_p \in \bm{S}_1, \\
\gamma_{iq\mid \bm{S}\setminus\{X_q\}} &= -\sum_{X_p\in\bm{S}_1}\tilde{\gamma}_{ip\mid \bm{S}_1\setminus\{X_p\}}\tau_{pq\mid do(\bm{S}_2 \setminus \{X_q\})}, \quad X_q \in \bm{S}_2.
\end{align*}
The first condition of the selective-door criterion implies
$\tau_{pq \mid do(\bm{S}_2 \setminus {X_q})} = 0$, and the second one implies $X_j \in \bm{S}_2$, which completes the proof.
\end{proof}

\begin{proof}[Proof of Corollaries \ref{col:DSA}, \ref{col:SDA}, and \ref{col:BDA}]
They follow immediately from Theorem~\ref{thm:selective-door criterion}.
\end{proof}

\begin{proof}[Proof of Proposition~\ref{prop:no confundings}]
Since it is clear that $(i)\Rightarrow(ii)\Rightarrow(iii)$ and $(ii)\Rightarrow(i)$, it suffices to prove $(iii)\Rightarrow(ii)$.
Fix a vertex $X_j \in \bm{S}$.
If there exists an opened back-door path from $X_j$ to $X_i$ that contains one or more v-structures, then there exists a back-door path $X_k \dashleftarrow X_i$ that contains no v-structures, where $X_k \in \bm{S}\setminus\{X_j\}$ is either the collider closest to $X_i$ along the original back-door path or a descendant of that collider.
This contradicts statement $(iii)$.
Hence, for any $X_j \in \bm{S}$, the set $\bm{S}\setminus\{X_j\}$ blocks every back-door path from $X_j$ to $X_i$.
Therefore, $(iii)\Rightarrow(ii)$ holds.
\end{proof}

\begin{proof}[Proof of Theorem~\ref{thm:PTB}]
See the proof of Theorem~\ref{thm:selective-door criterion}.
\end{proof}

\begin{proof}[Proof of Theorem~\ref{thm:Necessity almost everywhere}]
Here, we follow the notation used in the proof of Theorem~\ref{thm:selective-door criterion}.
Since the sufficiency follows from Theorem~\ref{thm:selective-door criterion}, we prove the necessity for almost all parameter points.
In the proof of Theorem~\ref{thm:selective-door criterion}, the following representation of biases is proved.
\begin{align*}
\gamma_{ip\mid \bm{S}\setminus\{X_p\}} &= \tilde{\gamma}_{ip\mid \bm{S}_1\setminus\{X_p\}}, \quad X_p \in \bm{S}_1, \\
\gamma_{iq\mid \bm{S}\setminus\{X_q\}} &= -\sum_{X_p\in\bm{S}_1}\tilde{\gamma}_{ip\mid \bm{S}_1\setminus\{X_p\}}\tau_{pq\mid do(\bm{S}_2 \setminus \{X_q\})}, \quad X_q \in \bm{S}_2.
\end{align*}
Since any subset of the parameter space defined by nontrivial algebraic equation constraints has Lebesgue measure zero, it suffices to show that these biases do not vanish identically as functions of the parameters.
In particular, if each $\tau_{pq \mid do(\bm{S}_2 \setminus \{X_q\})}$ vanishes identically, then the first condition of the selective-door criterion has to be satisfied.
Therefore, it is sufficient to show that each $\tilde{\gamma}_{ip\mid \bm{S}_1\setminus\{X_p\}}$ does not vanish identically.

Let $\tilde{\bm{X}}_{\bm{S}_1} = (\tilde{X}_{p_1}, \ldots, \tilde{X}_{p_M})^\top$ be the random vector consisting of the variables $\tilde{X}_p$, where the ordering is arbitrary.
Because an intercept term is included, the partial regression coefficients of interest (i.e., biases from the controlled total effects) can be expressed as
\begin{align*}
    \begin{pmatrix}
        \tilde{\gamma}_{ip_1\mid \bm{S}_1\setminus\{X_{p_1}\}} & \cdots & \tilde{\gamma}_{ip_M\mid \bm{S}_1\setminus\{X_{p_M}\}}    
    \end{pmatrix}^{\top}
    =\Var(\tilde{\bm{X}}_{\bm{S}_1})^{-1}\Cov(\tilde{\bm{X}}_{\bm{S}_1},\tilde{X}_i).
\end{align*}
Thus, by Cramer's rule, we have
\begin{align*}
    \tilde{\gamma}_{ip_1\mid \bm{S}_1\setminus\{X_{p_1}\}} = \frac{\det \bm{\Phi}_{p_1}}{\det \Var(\tilde{\bm{X}}_{\bm{S}_1})},
\end{align*}
where
\begin{align*}
    \bm{\Phi}_{p_1} =
    \begin{pmatrix}
        \Cov(\tilde{\bm{X}}_{\bm{S}_1},\tilde{X}_i) & \Cov(\tilde{\bm{X}}_{\bm{S}_1},\tilde{X}_{p_2}) & \Cov(\tilde{\bm{X}}_{\bm{S}_1},\tilde{X}_{p_3}) & \cdots & \Cov(\tilde{\bm{X}}_{\bm{S}_1},\tilde{X}_{p_M})
    \end{pmatrix}.
\end{align*}
Furthermore, by expanding the determinant along the first column, the numerator can be expressed as
\begin{align*}
    \det \bm{\Phi}_{p_1} = \sum_{k=1}^{M} (-1)^{k+1}\Cov(\tilde{X}_{p_k},\tilde{X}_i)\det \bm{\Phi}_{p_1}^{[k,1]},
\end{align*}
where $\bm{\Phi}_{p_1}^{[k,1]}$ denotes the submatrix of $\bm{\Phi}_{p_1}$ obtained by removing the $k$-th row and the first column.

We verify the following two statements.
\begin{enumerate}
    \item[(a)] $\Cov(\tilde{\bm{X}}_{\bm{S}_1},\tilde{X}_i)$ is not identically equal to the zero vector.
    \item[(b)] The minor $\det \bm{\Phi}^{[m,1]}_{p_1}$ corresponding to a component $\Cov(\tilde{X}_{p_m}, \tilde{X}_i)$ that does not vanish identically also does not vanish identically.
\end{enumerate}

(a). Note that the following equations hold from Proposition~\ref{prop:ancestral expansion}.
\begin{align*}
    \tilde{X}_i &=\sum_{X_l \in \bm{AN}(X_i)\setminus\bm{S}} \tau_{ik\mid do(\bm{S})}u_l+u_i, \\
    \tilde{X}_p &=\sum_{X_l \in \bm{AN}(X_p)\setminus\bm{S}_2} \tau_{pq\mid do(\bm{S}_2)}u_l+u_p.
\end{align*}
This implies that covariances among $\tilde{X}_i, \tilde{X}_{p_1}, \ldots, \tilde{X}_{p_M}$ are polynomial functions of the parameters.
From the definition of $\bm{S}_1$ and path-blocking, for a vertex $X_{p_m} \in \bm{S}_1$, there exists a back-door path $X_{p_m} \dashleftarrow X_i$ opened by both $\bm{S}\setminus \{X_{p_m}\}$ and the empty set.
Hence, at least one component of $\Cov(\tilde{\bm{X}}_{\bm{S}_1},\tilde{X}_i)$ does not vanish identically.

(b). If the first component $\Cov(\tilde{X}_{p_1},\tilde{X}_i)$ does not vanish identically, then, the statement holds immediately since $\det \bm{\Phi}^{[1,1]}_{p_1}$ is a principal minor of $\Var(\tilde{\bm{X}}_{\bm{S}_1})$.
Thus, we consider the case where $\Cov(\tilde{X}_{p_1}, \tilde{X}_i) \equiv 0$, that is, where each back-door path from $X_{p_1}$ to $X_i$ that is opened by $\bm{S}_1 \setminus \{X_{p_1}\}$ contains a collider.
Choose one such back-door path and consider the following sequence of vertices:
\[
    X_{p_1} \leftarrow\leftrightarrow\rightarrow X_{p_2}
    \leftarrow\leftrightarrow\rightarrow \cdots
    \leftarrow\leftrightarrow\rightarrow X_{p_m}
    \dashleftarrow X_i,
\]
where the back-door path $X_{p_m} \dashleftarrow X_i$ contains no v-structures, and each expression
$X_{k_1} \leftarrow\leftrightarrow\rightarrow X_{k_2}$ denotes a concatenated path consisting of
two directed paths sharing a common starting vertex, or a concatenation of two directed paths
connected by a bidirected edge.
Moreover, each of $X_{p_2}, \ldots, X_{p_m} \in \bm{S}_1$ is either a collider on the chosen back-door path
or a descendant of such a collider.
Finally, except for the vertices $X_{p_1}, \ldots, X_{p_m}$, none of the vertices on these paths
belongs to $\bm{S}$.

Then, $\Cov(\tilde{X}_{p_m}, \tilde{X}_i)$ does not vanish identically, and the corresponding  submatrix $\bm{\Phi}^{[m,1]}_{p_1}$ can be expressed as follows:
\begin{align*}
    \bm{\Phi}^{[m,1]}_{p_1} = 
    \begin{pmatrix}
        \tilde{\phi}_{p_1p_2} & \cdots & \tilde{\phi}_{p_1p_m} & \tilde{\phi}_{p_1p_{m+1}} & \cdots & \tilde{\phi}_{p_1p_M} \\
        \vdots & \ddots & \vdots & \vdots & \ddots & \vdots  \\
        \tilde{\phi}_{p_{m-1}p_2} & \cdots & \tilde{\phi}_{p_{m-1}p_m} & \tilde{\phi}_{p_{m-1}p_{m+1}} & \cdots & \tilde{\phi}_{p_{m-1}p_M} \\
        \tilde{\phi}_{p_{m+1}p_2} & \cdots & \tilde{\phi}_{p_{m+1}p_m} & \tilde{\phi}_{p_{m+1}p_{m+1}} & \cdots & \tilde{\phi}_{p_{m+1}p_M} \\
        \vdots & \ddots & \vdots & \vdots & \ddots & \vdots  \\
        \tilde{\phi}_{p_Mp_2} & \cdots & \tilde{\phi}_{p_Mp_m} & \tilde{\phi}_{p_Mp_{m+1}} & \cdots & \tilde{\phi}_{p_Mp_M}
    \end{pmatrix},
\end{align*}
where $\tilde{\phi}_{p_s p_t}\coloneq\Cov(\tilde{X}_{p_s}\tilde{X}_{p_t})$.
Due to the existence of the sequence of vertices considered above, it is clear that the determinant necessarily contains the product of diagonal entries that do not vanish identically.
For $1 < k < m-1$, although the entry $\tilde{\phi}_{p_k p_{k+1}} = \tilde{\phi}_{p_{k+1} p_k}$ is included in $\bm{\Phi}^{[m,1]}_{p_1}$, 
replacing it with the corresponding diagonal entry would cause a conflict with the indices of the other diagonal entries, 
so no other term in the determinant can identically cancel the product of the diagonal entries.
This means
$\det \bm{\Phi}^{[m,1]}_{p_1} \not \equiv 0$.

From facts (a) and (b), for almost all parameter points, the determinant of interest can be expressed as the inner product of two nonzero vectors with partially overlapping supports:
\begin{align*}
    \det \bm{\Phi}_{p_1}
    = \Cov(\tilde{\bm{X}}_{\bm{S}_1}, \tilde{X}_i)^\top \bm{g}(\Var(\tilde{\bm{X}}_{\bm{S}_1})),
\end{align*}
where $\bm{g}$ is a vector-valued polynomial function of the entries of $\Var(\tilde{\bm{X}}_{\bm{S}_1})$.
This implies that $\det \bm{\Phi}_{p_1}=0$ is a nontrivial algebraic equation constraint.
Therefore, $\tilde{\gamma}_{ip_1\mid \bm{S}_1\setminus\{X_{p_1}\}}
= \gamma_{ip_1\mid \bm{S}_1\setminus\{X_{p_1}\}}$
does not vanish for almost all parameter points.
Since the choice of the subscript $p_1$ was arbitrary, the proof is complete.
\end{proof}

\begin{proof}[Proof of Proposition \ref{prop:SDCvsRDC}]
We prove the claim by contradiction. 
Suppose, for the sake of contradiction, that a recanting district exists.
Then there exist two proper causal paths of the form
\begin{align*}
    X \rightarrow D_i \rightarrow \rightarrow Y, \\
    X \rightarrow D_j \rightarrow \rightarrow Y ,
\end{align*}
where $D_i$ and $D_j$ are connected by the bidirected edge $D_i \leftrightarrow D_j$, possibly with $D_i=D_j$. 
These two paths satisfy one of the following two cases:
\begin{enumerate}
    \item[(a)] $X \rightarrow D_i \rightarrow\rightarrow Y$ is blocked by $\bm{Z}$, but $X \rightarrow D_j \rightarrow\rightarrow Y$ is opened by $\bm{Z}$.
    \item[(b)] $X \rightarrow D_i \rightarrow\rightarrow Y$ is opened by $\bm{Z}$, but $X \rightarrow D_j \rightarrow\rightarrow Y$ is blocked by $\bm{Z}$.
\end{enumerate}
By symmetry of the structure, it suffices to consider only case (a).
Let $Z^*$ be the vertex in $\bm{Z}$ that is closest to $X$ along the path $X \rightarrow D_i \rightarrow\rightarrow Y$.  
Then, $Z^*$ is a post-treatment variable that violates the first condition of the selective-door criterion, contradicting the assumption. 
Therefore, the proposition follows.
\end{proof}

\end{document}